\documentclass[12pt,dvips]{amsart}
\usepackage{amsfonts, amssymb, latexsym, epsfig}
\usepackage{euler, amsfonts, amssymb, latexsym, epsfig}
\usepackage{euler, amsfonts, amssymb, latexsym, epsfig,epic}

\newtheorem{Theorem}{Theorem}[section] 

\newtheorem{Proposition}[Theorem]{Proposition} 
\newtheorem{lemma}[Theorem]{Lemma}

\newtheorem{corollary}[Theorem]{Corollary}
\newtheorem{problem}[Theorem]{Problem}
\newtheorem{conjecture}[Theorem]{Conjecture}              
\newtheorem*{Corollary*}{Corollary}
 
\newtheorem*{Theorem*}{Theorem}

\theoremstyle{remark}
\newtheorem{Example}[Theorem]{Example}

\newcommand\Spec{\operatorname{Spec}}
\newcommand\zvect{{\bf z}}

\newcommand{\idperm}{{\rm id}}
\newcommand{\Hom}{\operatorname{Hom}}
\newcommand{\Ext}{\operatorname{Ext}}
\newcommand{\Tor}{\operatorname{Tor}}
\newcommand{\codim}{\operatorname{codim}}

\newcommand{\ourfield}{\mathbb{C}}

\theoremstyle{plain}

\begin{document}
\pagestyle{plain}

\title{Governing Singularities of Schubert varieties}
\author{Alexander Woo}
\address{Department of Mathematics, Mathematical Sciences Building, One Shields Ave., 
University of California, Davis, CA, 95616, USA}
\email{awoo@math.ucdavis.edu}

\author{Alexander Yong}
\address{Department of Mathematics, University of Minnesota, 
Minneapolis, MN 55455, USA; Department of Statistics {\it and}
the Fields Institute, University of Toronto, Toronto, Ontario, M5T 3J1,
Canada}

\email{ayong@math.umn.edu, ayong@fields.utoronto.ca}

\date{June 29, 2006}
\begin{abstract}
We present a combinatorial and computational commutative algebra
methodology for studying singularities of Schubert varieties of 
flag manifolds. 

We define 
the combinatorial notion of \emph{interval pattern avoidance}. 
For ``reasonable'' invariants ${\mathcal P}$ 
of singularities, we geometrically prove 
that this governs (1) the 
${\mathcal P}$-locus of a Schubert variety, and (2) which
Schubert varieties are globally not ${\mathcal P}$. The prototypical case
is ${\mathcal P}=$``singular''; \emph{classical} pattern avoidance
applies admirably for this choice [Lakshmibai-Sandhya'90], but 
is insufficient in general. 

Our approach is analyzed for some common invariants, including
Kazhdan-Lusztig polynomials, multiplicity, factoriality, and
Gorensteinness, extending [Woo-Yong'04]; the description of the singular locus
(which was independently proved by [Billey-Warrington '03], [Cortez
'03], [Kassel-Lascoux-Reutenauer'03], [Manivel'01]) is also thus
reinterpreted.

Our methods are amenable to computer experimentation, based on
computing with \emph{Kazhdan-Lusztig ideals} (a class of generalized
determinantal ideals) using {\tt Macaulay 2}. This feature is supplemented
by a collection of open problems and conjectures.
\end{abstract}
\maketitle
\vspace{-.5in}
\tableofcontents
\section{Overview}

Let $X_w$ be the Schubert variety of the complete flag variety ${\rm
Flags}({\mathbb C}^n)$ associated to a permutation~$w$ in the
symmetric group~$S_n$.  One would like to describe and classify the
singularities of $X_w$, as well as calculate invariants measuring
their complexity.  Solutions to such problems typically
require techniques from and have important applications to geometry,
representation theory, and associated combinatorics.  Two recent
surveys of some work in this area are~\cite{BL00, brion:book}.

In this paper, we formulate a new combinatorial notion, a
generalization of pattern avoidance we call \emph{interval pattern
avoidance}; we then use this idea to explore the singularities of
Schubert varieties and their local invariants.  The well-known
Kazhdan-Lusztig polynomials show up as one local invariant, since
their coefficients are the Betti numbers for the local intersection
cohomology of the singularities.  Indeed, a desire to further
understand the combinatorics of Kazhdan-Lusztig polynomials is one
source of motivation (and application) for this present work.
However, there are many other noteworthy invariants of singularities,
including factoriality, multiplicity, Gorensteinness, and
Cohen-Macaulay type.  We provide a uniform language to study such
\emph{semicontinuously stable} invariants, in an attempt to gain
further insight into the singularities of Schubert varieties.

Informally, our principal thesis is that, for any of these
``reasonable'' local invariants of singularities of Schubert
varieties, the question of where it assumes a particular value has a
natural answer in terms of interval pattern avoidance.  Our main
result (Theorem~\ref{thm:poset}) is a precise version of this
assertion, together with a geometric explanation; proofs are given in
Section~4.

The two most basic
problems about singularities of specific Schubert varieties are
\begin{itemize}
\item Which $X_w$ are singular?
\item Where is $X_w$ singular?
\end{itemize}

These questions have been answered.  Following upon a geometric
characterization by Ryan~\cite{Ryan} and an earlier combinatorial
characterization by Wolper~\cite{Wolper}, V.~Lakshmibai and B.~Sandhya
\cite{LS90} gave a simple characterization of singular Schubert
varieties in terms of the combinatorial notion of {\it pattern
avoidance}: $X_{w}$ is smooth if and only if $w$ avoids the patterns
$3412$ and $4231$; see the definitions in Section~2. They also
conjectured an explicit description for the singular locus of $X_w$ in
terms of pattern avoidance. This conjecture was solved independently
by several groups~\cite{billey.warrington, Cortez, KLR, manivel1}
around 2000.  We reinterpret this result in terms of interval pattern
avoidance.

Although much is known concerning general properties of singularities
of Schubert varieties, little more is known for properties which not
all Schubert varieties hold in common.  Thanks to fundamental work
during the 1980s including that of C. DeConcini and V.~Lakshmibai
~\cite{DeCon-Lak}, and S.~Ramanan and A.~Ramanathan~\cite{RR, Ram85},
we know that all Schubert varieties are Cohen-Macaulay and normal.  In
addition, A.~Cortez~\cite{Cortez} and L.~Manivel~\cite{manivel2}
independently described the neighborhoods of \emph{generic} points in
the singular locus of a Schubert variety; understanding where and how
these neighborhoods change at special points of the singular locus is
a core theme in our present investigations.  More recently
in~\cite{WY:Gor}, we determined which Schubert varieties are
Gorenstein; we introduced a notion there called Bruhat-restricted
pattern avoidance, and interval pattern avoidance is a further
generalization which has the advantage of a geometric interpretation.
We further pursue below the question of where a non-Gorenstein
Schubert variety is Gorenstein, along with analogous questions for
other local properties.

Analysis of specific questions from this viewpoint suggests new algebraic,
geometric, and combinatorial questions and conjectures which we explore
computationally using {\tt Macaulay~2}~\cite{macaulay}. This is
explained in Sections~5 and~6.  The associated commutative algebra is
that of \emph{Kazhdan-Lusztig ideals} (a class of ideals generalizing
classical determinantal ideals); this commutative algebra is explicated in Section~3.

This report was written in part to help facilitate activities at the
``Workshop on combinatorial and computational commutative algebra'' 
(Fields Institute, July-August 2006). The workshop advances the use of 
computer algebra systems such as {\tt Macaulay~2}. 
We wrote the {\tt Macaulay~2} code 
{\tt Schubsingular} as an exploratory complement to this 
paper.\footnote{Available at the authors' websites.} 

For simplicity, this paper focuses on the complete flag manifold in
type~$A$. This allows us to emphasize links to the traditional study
of determinantal ideals in commutative algebra and avoid the need for
terminology from the theory of algebraic groups.  However, the ideas
below can be extended with appropriate modifications to the other root
systems and partial flag manifolds.  Finally,
although in this paper we work over $\mathbb{C}$ for convenience, our
results are valid over any field $\Bbbk$ of any characteristic except
as noted.

\section{The main definitions and theorem}

\subsection{Interval pattern avoidance}

Let $v\in S_m$ and $w\in S_n$ be two permutations, where $m\leq n$.
We say $v$ {\bf embeds in} $w$ if there exist indices
$1\leq\phi_1<\phi_2<\ldots<\phi_m\leq n$ such that $w(\phi_1),
w(\phi_2), \ldots, w(\phi_m)$ are in the same relative order as
$v(1),\ldots, v(m)$.  In other words, we require that $w(\phi_j)<w(\phi_k)$
if and only if $v(j)<v(k)$.  The permutation $w$ is said to
(classically) {\bf avoid} $v$ if no such embedding exists.

Recall that {\bf Bruhat order}, which we denote by $\leq$, is the
partial order on $S_m$ defined by declaring that $u\leq v$ if
$v=u(i\leftrightarrow j)$ and $\ell(v)>\ell(u)$, and taking the
reflexive transitive closure.  Here, $(i\leftrightarrow j)$ is the
transposition switching positions $i$ and $j$, and $\ell(v)$ denotes
the {\bf Coxeter length} of $v$, which is the length of any reduced
expression for $v$ as a product of simple reflections
$s_i=(i\leftrightarrow i+1)$.  Alternatively, $\ell(v)$ is also the
number of {\bf inversions} of $v$; inversions are pairs $i,j$ with
$1\leq i<j\leq m$ such that $v(i)>v(j)$.  Bruhat order is a partial
order graded by Coxeter length.

We now give our main definition.  Let $[u,v]$ and $[x,w]$ be intervals
in the Bruhat orders on $S_m$ and $S_n$ respectively.  We say that
$[u,v]$ {\bf (interval) pattern embeds in} $[x,w]$ if there is a
common embedding $\Phi=(\phi_1,\ldots,\phi_m)$ of $u$ into $x$ and $v$
into $w$, where the entries of $x$ and $w$ outside of $\Phi$ agree,
and, furthermore, $[u,v]$ and $[x,w]$ are isomorphic as posets.

Note that the first two requirements already determine $x$ given $u$,
$v$, $w$, and $\Phi$.  To be precise, for a permutation $\sigma\in
S_m$, let $\Phi(\sigma)\in S_n$ be the permutation where
$\Phi(\sigma)(\phi_j)=w(\phi_{(v^{-1}\sigma)(j)})$, and
$\Phi(\sigma)(k)=w(k)$ if $k\neq\phi_j$ for $1\leq j\leq m$.  Then the
first two requirements force $x$ to be equal to $\Phi(\sigma)$.
Therefore, for convenience, we sometimes drop $x=\Phi(u)$ and say
that $[u,v]$ embeds in $w$ if $[u,v]$ embeds in $[\Phi(u),w]$.  We
also say simply that $w$ {\bf (interval) (pattern) avoids}
$[u,v]$ if there are no interval pattern embeddings of $[u,v]$ into
$[x,w]$ for any $x\leq w$.

The following lemma gives a simple criterion for checking if a pattern
embedding actually produces an interval pattern embedding.  Its proof
is simple and we omit it.

\begin{lemma}
\label{lemma:embed_iso}
An embedding $\Phi$ of $[u,v]$ into $[\Phi(u),w]$ is an interval
pattern embedding if and only if
$\ell(v)-\ell(u)=\ell(w)-\ell(\Phi(u))$.
\end{lemma}

\begin{Example}
\label{exa:35142}
Let $v=35142= s_2 s_1 s_4 s_3 s_2 s_4$ and $u = 13524 = s_2 s_4 s_3$.
Note $u\leq v$, and $\ell(v)-\ell(u)=3$.  Now let $\Phi$ be the
embedding of $v$ into $w=\underline{5} 8 9 \underline{7} \underline{1}
\underline{6} 2 \underline{3} 4$ where the underlined positions
indicate the embedding, which in symbols is given by $\phi_1 =1,
\phi_2 = 4, \phi_3 =5, \phi_4 =6, \phi_5 =8$. Then $\Phi(u) =
\underline{1} 8 9 \underline{5} \underline{7} \underline{3} 2
\underline{6} 4$. The reader can check that $\ell(w) = 24$ and
$\ell(\Phi(u))=21$. Therefore this is an embedding of $[u,v]$ into
$[\Phi(u),w]$.\qed
\end{Example}

\begin{Example}
\label{exa:2413}
Let $v=2413= s_1 s_3 s_2$, $u=2143 = s_1 s_3$, and note that
$\ell(v)-\ell(u)=1$.  Let $w = 265314$; note there are two
embeddings $\Phi_1$ and $\Phi_2$ of $v$ into $w$, represented
respectively by the underlinings $\underline{2} \underline{6} 5
3\underline{1} \underline{4}$ and $\underline{2} 6 \underline{5}
3\underline{1}\underline{4}$.

Neither of these embeddings induce an embedding of $[u,v]$.  We have that
$\Phi_1(u) = 215364$ and $\ell(w)-\ell(\Phi_1(u))=5\neq
\ell(v)-\ell(u)=1$.  (Note that $\ell(w)=9$.)  For $\Phi_2$,
$\Phi_2(u) = 261354$ and $\ell(w)-\ell(\Phi_2(u))=3$ which again
differs from $\ell(v)-\ell(u)=1$. Hence $w$ in fact
interval pattern avoids $[u,v]$, even though it does not classically
pattern avoid $v$.\qed
\end{Example}

Two further lemmas follow immediately from our definition.

\begin{lemma}
If $\Phi$ gives an embedding of $[u,v]$ into $[\Phi(u),w]$, then
$\Phi$ also gives an embedding of $[u^\prime,v]$ into
$[\Phi(u^\prime), w]$ for any $u^\prime$ such that $u\leq u^\prime\leq
v$.
\end{lemma}

\begin{lemma}
If $w$ avoids $[u,v]$, then $w$ avoids $[u^\prime,v]$ for any
$u^\prime\leq u$.
\end{lemma}

\subsection{Semicontinuously stable properties}

We are interested in local properties of Schubert varieties that are
semicontinuous, meaning that they hold on closed subsets of any
variety, and are preserved under products with affine space. We call such
a property ${\mathcal P}$ {\bf semicontinuously stable}.  For example,
\[\{\mbox{semicontinuously stable }{\mathcal P}\} =
\left\{
\begin{array}{c}
\mbox{singular, non-Gorenstein,} \\ 
\mbox{non-factorial, dimension of $i$-th local intersection}\\ 
\mbox{cohomology group $\geq k$, Cohen-Macaulay type $\geq k$,}\\  
\mbox{multiplicity $\geq k$},\ldots
\end{array}
\right\}.
\]
For us, ``reasonable'' invariants of singularities are properties such that
they, or their negations, are semicontinuously stable. (Actually, at present there 
is no general result that the property 
${\mathcal P}=$``dimension of $i$-th local intersection cohomology group 
$\geq k$'' is semicontinuously stable, but for Schubert varieties 
this is known to be true~\cite{Irving}.)

We desire a common combinatorial language to describe the {\bf
${\mathcal P}$-locus}, the closed subset of a Schubert variety $X_w$
at which the local property ${\mathcal P}$ holds. (As explained in
Section~3 below, the $\mathcal{P}$-locus is a union of Schubert
subvarieties, and it suffices to consider this question at the
$T$-fixed points $e_x$ for $x\leq w$.)  As it turns out, classical
pattern avoidance is insufficient in general for this purpose. We
first observed this for Gorensteinness in~\cite{WY:Gor}.  (See
Example~\ref{exa:Gor_not_pattern} below.)  It was also there that we
first noticed the phenomenon of Bruhat-restricted/interval pattern avoidance,
which suggested the present study.


To connect the combinatorics of interval pattern avoidance to
the geometry of Schubert varieties, we need a little more notation and
terminology.  Consider the set
\[\mathfrak{S}=\{[u,v]: u\leq v \mbox{ in some $S_r$}\}\subseteq S_\infty \times S_{\infty}\]
where $S_\infty  = \bigcup_{r\geq 1} S_r$.

Define $\prec$ to be the partial order on $\mathfrak{S}$ generated by
the two types of relations
\begin{enumerate}
\item $[u,v]\prec [x,w]$ if there is an interval pattern embedding of
$[u,v]$ into $[x,w]$, and
\item $[u,v]\prec [u^\prime, v]$ if $u^\prime\leq u$.
\end{enumerate} 
 
An {\bf upper order ideal} ${\mathcal I}$ (under the partial order
$\prec$) is a subset of $\mathfrak{S}$ such that, if $[u,v]\in
{\mathcal I}$ and $[u,v]\prec [x,w]$, then $[x,w]\in {\mathcal I}$.

We are now ready to state the precise version of our main idea from
Section~1:
\begin{Theorem}
\label{thm:poset}
Let ${\mathcal P}$ be a semicontinuously stable property. The set of
intervals $\{[u,v]\}\subseteq \mathfrak{S}$ such that ${\mathcal P}$
holds at the $T$-fixed point $e_u$ on the Schubert variety $X_v$ is an
upper order ideal ${\mathcal I}_{\mathcal P}$ under $\prec$.
\end{Theorem}

We also wish to characterize Schubert varieties that globally avoid
${\mathcal P}$, or, in other words, those Schubert varieties for which
${\mathcal P}$ does not hold at any point, in analogy with the
theorems for smoothness~\cite{LS90} and Gorensteinness~\cite{WY:Gor}.
The following corollary says that this can be done in terms of
interval pattern avoidance.

\begin{corollary}
\label{cor:brpa}
Let $\mathcal{P}$ be a semicontinuously stable property. Then the set
of permutations $w$ such that $\mathcal{P}$ does not hold at any point
of $X_w$ is the set of permutations $w$ that avoid all the intervals
$[u_i,v_i]$ constituting some (possibly infinite) set
$A_\mathcal{P}\subseteq \mathfrak{S}$.
\end{corollary}

The corollary is false in general for classical pattern avoidance, as
the following example illustrates:

\begin{Example}
\label{exa:Gor_not_pattern}
The Schubert variety $X_{42513}\subseteq {\rm Flags}({\mathbb C}^5)$
is not Gorenstein (see Theorem~\ref{thm:Gor_char}). However
$X_{526413}\subseteq {\rm Flags}({\mathbb C}^6)$ is Gorenstein
even though $42513$ embeds into $\underline{526}4\underline{13}$
at the underlined positions. So Gorensteinness cannot be characterized
using classical pattern avoidance.\qed
\end{Example}

We speculate that for any semicontinuously stable property ${\mathcal
P}$, $\mathcal{I}_\mathcal{P}$ and $A_\mathcal{P}$ respectively
provide natural answers to the problems of
\begin{itemize}
\item determining the ${\mathcal P}$-locus 
of $X_w$ and
\item characterizing which Schubert varieties $X_w$ 
globally avoid ${\mathcal P}$.
\end{itemize}
Therefore, we expect interval pattern avoidance to be a
useful framework for studying these questions, both in principle, as
established by the above results, and practice, as evidenced by the
examples in Section~6.

In Sections~3 and~4, we introduce Kazhdan-Lusztig ideals, using them
to explain and prove Theorem~\ref{thm:poset} and its corollary. After
this we will proceed to describe $\mathcal{I}_\mathcal{P}$ and
$A_\mathcal{P}$ for the properties $\mathcal{P}$ for which they are
known, and explain how $\mathcal{I}_\mathcal{P}$ and $A_\mathcal{P}$
might be computed for some other properties.  Note that
$\mathcal{I}_\mathcal{P}$ and $A_\mathcal{P}$ may vary depending on
the base field $\Bbbk$ for some properties $\mathcal{P}$.

\section{Schubert varieties and Kazhdan-Lusztig ideals}

\subsection{Schubert definitions} 
Let $G=GL_{n}({\ourfield})$ denote the group of
invertible $n\times n$ matrices with entries in $\ourfield$, and let $B,
B_{-}, T\subseteq G$ denote the subgroups of upper triangular, lower
triangular and diagonal matrices respectively.  The {\bf complete flag
variety} is ${\rm Flags}({\ourfield}^n)=G/B$; upon choosing a basis of
$\ourfield^n$, a point $gB\in G/B$ is naturally identified with a complete
flag $F_{\bullet}: \langle 0\rangle\subsetneq F_1\subsetneq
F_2\subsetneq\cdots\subsetneq F_n ={\ourfield}^n$ by allowing $F_i$ to be
the span of the first $i$ columns of any coset representative
of $gB$.

The flag variety has a {\bf Bruhat decomposition}
$$G/B=\coprod_{w\in S_n} BwB/B,$$
where we think of $w$ as the permutation matrix with $1$'s at $(w(i),i)$
and $0$'s elsewhere.  The Zariski closure of the {\bf Schubert cell} 
$X_{w}^{\circ}:=BwB/B$ is the Schubert variety 
$X_{w}:={\overline{X_{w}^{\circ}}}$. The
Schubert cell $X^\circ_w$ is isomorphic to affine space~${\mathbb A}^{\ell(w)}$.  Moreover, 
each Schubert variety is the disjoint union of Schubert cells
$$X_{w}=\coprod_{x\leq w} X_{x}^{\circ}.$$ 
Our conventions have been chosen so that the dimension of $X_w$ is
$\ell(w)$. In particular, 
$X_{\mathrm{id}}$ is a point, and $X_{w_0}=G/B$, where
$w_0$ denotes the permutation such that $w_0(i)=n+1-i$.

The $T$-fixed points of $X_w$ (under the left action of $T$ on $G/B$) 
are $e_x:=xB/B$
for $x\leq w$; these are known as {\bf Schubert points} and represent
the flags corresponding to permutation matrices. Every point on
a Schubert variety is in the $B$-orbit of some Schubert point, and the
$B$-action gives an isomorphism between a local neighborhood of any
point with a local neighborhood of a Schubert point.  
Therefore, in studying local questions about singularities of Schubert
varieties, we may restrict attention to the Schubert points.

\subsection{Affine neighborhoods, explicitly}

Let $M_n$ be the space of $n\times n$ matrices over $\ourfield$; we
think of $M_n$ as a variety with coordinate ring $\ourfield[{\bf z}]$
where ${\bf z}:=\{z_{n-i+1,j}\}_{i,j=1}^{n}$ are the entries of a
generic matrix~$Z$. (Note we are labeling our variables so that
$z_{11}$ is at the southwest corner (at row~$n$ and column~$1$) of
the generic matrix.)

For a permutation $x\in S_n$, let $Z^{(x)}$ be the generic matrix
obtained when we specialize $Z$ by setting $z_{n-x(i)+1,i}=1$ for all
$i$, and $z_{n-x(i)+1,a}=0$ and $z_{n-b+1,i}=0$ for $a>i$ and
$b<x(i)$. Let ${\bf z}^{(x)}\subseteq {\bf z}$ denote the other
(unspecialized) variables. Let $Z^{(x)}_{ij}$ denote the
southwest submatrix of $Z^{(x)}$ with northeast corner $(i,j)$; this
matrix has $n-i+1$ rows and $j$ columns.  Furthermore, let
$R^w=\left[r_{ij}^w\right]_{i,j=1}^n$ be the {\bf rank matrix}, in
which each $r_{ij}^w$ equals the number of $1$'s to the southwest of
$(i,j)$ in~$w$:
$$r_{ij}^w=\#\{k\mid w(k)\geq i \mbox{ and } k\leq j\}.$$ Let the
{\bf Kazhdan-Lusztig ideal} $I_{x,w}$ be generated by the size 
$1+r_{ij}^w$ minors of $Z_{ij}^{(x)}$, over all
possible $i$ and $j$. Let
\[{\mathcal N}_{x,w}:=\Spec(\ourfield[{\bf z}^{(x)}]/I_{x,w})\]
be the associated affine scheme. (See Example~\ref{exa:35142var} below.)

We formulate our proof of Theorem~\ref{thm:poset} and the computations
in Sections~5--6 using the following fact. (Experts will find
the ideas contained herein familiar.)

\begin{Proposition}
${\mathcal N}_{x,w} \times \mathbb{A}^{\ell(x)}$ is isomorphic to an
affine neighborhood of $X_{w}$ at $e_x$. In particular, if ${\mathcal
P}$ is a semicontinuously stable property, then $\mathcal{P}$ holds at
$e_x$ on $X_w$ if and only if $\mathcal{P}$ holds at the origin
$\mathbf{0}$ on ${\mathcal N}_{x,w}$.
\end{Proposition}

\noindent
\emph{Proof and discussion:}
An affine neighborhood of $e_x$ in the
flag variety is given by $x\Omega^\circ_\idperm$, where, in general,
$\Omega^\circ_u$ is the {\bf opposite Schubert cell} defined by
$\Omega^\circ_u:=B_{-}uB/B=w_0X^\circ_{w_0u}$.  To study $X_w$ locally at $e_x$,
we therefore need only understand $X_w\cap x\Omega^\circ_\idperm$.  

We now proceed to describe explicit coordinates for $x\Omega^\circ_\idperm$ and
equations for $X_w\cap x\Omega^\circ_\idperm$ in terms of these coordinates.

Let $\pi: G\rightarrow G/B$ be the natural quotient map.  The map
$\pi$ has a local section $\sigma$ over $\Omega^\circ_\idperm$ with
$\sigma(F_\bullet)$ being the unique coset representative of
$F_\bullet$ which is unit lower triangular ($1$'s are on the main
diagonal).  The map $x\sigma x^{-1}$ is then a local section over
$x\Omega^\circ_\idperm$; therefore we have that $$X_w\cap
x\Omega^\circ_{\idperm}\cong\pi^{-1}(X_w)\cap x\sigma
x^{-1}(x\Omega^\circ_{\idperm}),$$ where the latter can be considered
as a subvariety of $M_n$.

The following lemma was first stated by D.~Kazhdan and G.~Lusztig
(whence our terminology for $I_{x,w}$). It holds for the flag
varieties of all algebraic groups.  For completeness, we give an
explicit description of the isomorphism in our $GL_{n}({\ourfield})$
case.

\begin{lemma}[Lemma~A.4 of \cite{KL}] 
\label{lemma:KL}
$X_w\cap x\Omega^\circ_\idperm\cong (X_w \cap \Omega^\circ_x)
\times \mathbb{A}^{\ell(x)}.$
\end{lemma}
\begin{proof}
The map $x\sigma x^{-1}$ sends $\Omega^\circ_x$ to the set of matrices
with $1$'s at $(x(i),i)$ for $1\leq i\leq n$, $0$'s to the right and
above these $1$'s, and arbitrary entries elsewhere.  Now identify
$\mathbb{A}^{\ell(x)}$ with the space of unit upper triangular
matrices $m=\left[m_{ij}\right]_{i,j=1}^n$ for which
$m_{ij}=0$ (for $i<j$) unless $x(i)>x(j)$.  It is easy to check that
the map $\eta:x\sigma x^{-1}(\Omega^\circ_x)\times
\mathbb{A}^{\ell(x)}\rightarrow x\sigma x^{-1}(x\Omega^\circ_\idperm)$
given by $\eta(a,m)=ma$ (where we have matrix multiplication on the
right hand side) is an isomorphism.  Now notice $\eta$ restricts as
desired to any Schubert variety $X_w$; if $a\in X_w \cap
\Omega^\circ_x$, then $\pi(\eta(x\sigma x^{-1}(a),m))\in X_w\cap
x\Omega^\circ_\idperm$, since $X_w$ is closed under the action of $B$.
\end{proof}

Let $\mathcal{N}'_{x,w}$ denote the variety $X_w\cap \Omega^\circ_x$.
In view of Lemma~\ref{lemma:KL}, that it 
suffices to study these varieties to 
understand semicontinuously stable properties of $X_w$. 
We want to show ${\mathcal N}'_{x,w}\cong {\mathcal N}_{x,w}$.

To do this, we want explicit equations in coordinates for
$\mathcal{N}'_{x,w}$.  Since $x\sigma x^{-1}:\Omega^\circ_x\rightarrow G$
is a section of the map $\pi$ and hence an injection, we have
$$\mathcal{N}'_{x,w}\cong x\sigma x^{-1}(\Omega^\circ_x)\cap \pi^{-1}(X_w).$$

One coordinate ring for $GL_{n}({\ourfield})$ is
$\ourfield[\zvect,\det^{-1}(\zvect)]$, where
$\zvect:=(z_{n-i+1,j})_{i,j=1}^n$ are the entries of a generic
\emph{invertible} matrix $Z$.  With these coordinates, the defining ideal for
$x\sigma x^{-1}(\Omega^\circ_x)$ is generated by the polynomials
$z_{n-x(i)+1,i}-1$ and monomials of the form $z_{n-x(i)+1,a}$ and
$z_{n-b+1,i}$ for $a>i$ and $b<x(i)$; we denote this ideal $J_x$.
Fulton~\cite{Fulton:Duke92} showed that the ideal $I_w$ defining
$\pi^{-1}(X_w)$ (scheme-theoretically) is generated by the size
$1+r_{ij}^w$ minors of $Z_{ij}$, over all possible $i$ and $j$; the
closure of $\pi^{-1}(X_w)$ in $M_n$ defined by $\ourfield[\mathbf{z}]/I_w$
is known as the {\bf matrix Schubert variety}.  Actually, as
Fulton~\cite{Fulton:Duke92} explains, $I_w$ is generated by a much
smaller subset of these minors, corresponding to the essential set
conditions; we describe this in the example below.\footnote{Our
conventions differ from those of \cite{Fulton:Duke92,grobGeom,
WY:Gor}; our equations define their (matrix) Schubert varieties for
$w_0w^{-1}$.}  Therefore,
$$\mathcal{N}'_{x,w}\cong\Spec(\ourfield[\zvect]/(I_w+J_x)).$$
(Note that $\det(\zvect)=1$ by $J_x$.)

In practice, to reduce the number of variables, instead of working in
a generic matrix $Z$, we first quotient by $J_x$ and work in the
generic matrix $Z^{(x)}$.  The image of $I_w$ in $\ourfield[{\bf
z}^{(x)}]$ under this quotient by $J_x$ is precisely $I_{x,w}$. Therefore
${\mathcal N}_{x,w}\cong {\mathcal N}'_{x,w}$ and the result follows.\qed 

The following is an immediate corollary.  However, it is far from
obvious if one looks only at the generators of $I_{x,w}$.

\begin{corollary}
${\mathcal N}_{x,w}$ is reduced and irreducible of dimension $\ell(w)-\ell(x)$.
\end{corollary}

The following gives an example of the theorem as well as its
proof and discussion above:

\begin{Example}
\label{exa:35142var}
Let $w=35142\in S_5$; then 
$$
Z=\left(\begin{matrix}
z_{51} & z_{52} & z_{53} & z_{54} & z_{55} \\
z_{41} & z_{42} & z_{43} & z_{44} & z_{45} \\
z_{31} & z_{32} & z_{33} & z_{34} & z_{35} \\
z_{21} & z_{22} & z_{23} & z_{24} & z_{25} \\
z_{11} & z_{12} & z_{13} & z_{14} & z_{15} 
\end{matrix}\right), \ 
w=\left(\begin{matrix}
0 & 0 & 1 & 0 & 0 \\
0 & 0 & 0 & 0 & 1 \\
1 & 0 & 0 & 0 & 0 \\
0 & 0 & 0 & 1 & 0 \\
0 & 1 & 0 & 0 & 0 
\end{matrix}\right) \mbox{ and } 
R_w = \left(\begin{matrix}
1 & 2 & 3 & 4 & 5 \\
1 & 2 & 2 & 3 & 4 \\
1 & 2 & 2 & 3 & 3 \\
0 & 1 & 1 & 2 & 2 \\
0 & 1 & 1 & 1 & 1 
\end{matrix}\right).
$$
Drawing ``hooks'' to the right and above every $1$ in $w$ defines the
{\bf diagram} of $w$, which is the set of positions not in any hook.
Here, the diagram is the set $\{(2,3), (4,1), (4,3), (5,1)\}$. The
{\bf essential set} consists of the northeast most boxes in each
connected component of the diagram, in this case, $(2,3)$, 
$(4,1)$ and $(4,3)$. Fulton~\cite{Fulton:Duke92} showed that the minors arising
from considering just these three positions generate all of $I_w$, so
\begin{multline}\nonumber
I_{w}= \left\langle z_{11}, z_{21},
\left|
\begin{matrix}
z_{31} & z_{32} & z_{33} \\
z_{21} & z_{22} & z_{23} \\
z_{11} & z_{12} & z_{13} 
\end{matrix}
\right|,
\left|
\begin{matrix}
z_{41} & z_{42} & z_{43} \\
z_{21} & z_{22} & z_{23} \\
z_{11} & z_{12} & z_{13} 
\end{matrix}
\right|,
\left|
\begin{matrix}
z_{41} & z_{42} & z_{43} \\
z_{31} & z_{32} & z_{33} \\
z_{11} & z_{12} & z_{13} 
\end{matrix}
\right|,
\left|
\begin{matrix}
z_{41} & z_{42} & z_{43} \\
z_{31} & z_{32} & z_{33} \\
z_{21} & z_{22} & z_{23} 
\end{matrix}
\right|,\right. \\
\left.\left|
\begin{matrix}
z_{21} & z_{22} \\
z_{11} & z_{12}
\end{matrix}
\right|,
\left|
\begin{matrix}
z_{21} & z_{23} \\
z_{11} & z_{13}
\end{matrix}
\right|,
\left|
\begin{matrix}
z_{22} & z_{23} \\
z_{12} & z_{13}
\end{matrix}
\right|\right\rangle.
\end{multline}
(The reader may find it helpful to argue why
all the other minors in $I_w$ are in the ideal generated by only these minors.)

Let $x=13254\leq w$; then a generic matrix of $x\sigma
x^{-1}(\Omega^\circ_x)$ is
$$
\left(\begin{matrix}
1      & 0      & 0      & 0      & 0 \\
z_{41} & 0      & 1      & 0      & 0 \\
z_{31} & 1      & 0      & 0      & 0 \\
z_{21} & z_{22} & z_{23} & 0      & 1 \\
z_{11} & z_{12} & z_{13} & 1      & 0 
\end{matrix}\right).
$$
We set $z_{51}=z_{32}=z_{43}=z_{14}=z_{25}=1$,
and all other variables \emph{except} 
$z_{11}, z_{12}, z_{13}, z_{21}, 
z_{22},$ $z_{23}, z_{31}, z_{41}$ to 0, resulting in the Kazhdan-Lusztig ideal
\begin{eqnarray}\nonumber
I_{x,w} & = & \langle z_{11},z_{21}, -z_{11} z_{23} + z_{21} z_{13} + 
z_{31} z_{12} z_{23} - z_{31} z_{13} z_{22},  z_{11} z_{22} - z_{21} z_{12} + z_{41} z_{12} z_{23} \\ \nonumber 
& & - z_{41} z_{13} z_{22}, z_{11} - z_{31} z_{12} - z_{41} z_{13}, 
z_{21} - z_{31} z_{22} - z_{41} z_{23}, z_{11} z_{22}-z_{21} z_{12}, \\ \nonumber 
& & z_{11} z_{23}-z_{21} z_{13}, z_{12} z_{23}-z_{22}z_{13}\rangle. 
\end{eqnarray}\qed
\end{Example}

We remark that any matrix Schubert variety can be realized
(up to a product with affine space) as a particular
${\mathcal N}_{x,w}$. Specifically, for $w\in S_n$, consider its natural
embedding into $S_{2n}$ and let $x=id\in S_{2n}$. Then it is not hard
to check that the generators of $I_{\idperm,w}$ are the same as the
generators of $I_w$, except that the ideal $I_w$ is defined in a ring
with more variables.

\section{A local isomorphism and the proof of Theorem \ref{thm:poset}, Corollary~\ref{cor:brpa}}

The following lemma explains the role of semicontinuity
in our discussion.

\begin{lemma}
\label{lemma:semicont}
If $u^\prime\leq u$ and $X_v$ has a semicontinuously stable property
$\mathcal{P}$ at $e_u$, then $X_v$ must also have this property at
$e_{u^\prime}$.  In particular, if property $\mathcal{P}$ fails to 
hold at $e_{\idperm}$, it will not
hold for any point on $X_v$.
\end{lemma}

\begin{proof}
It suffices to consider the case where $u$ actually covers $u^\prime$
in Bruhat order.  There is a $\mathbb{P}^1$ connecting $e_{u^\prime}$
and $e_u$.  In this $\mathbb{P}^1$, the generic point is in the
$B$-orbit of $e_u$, and hence has a neighborhood isomorphic to a
neighborhood of $e_u$.  The one special point is $e_u^\prime$.  Since
$\mathcal{P}$ is semicontinuous and holds on a generic point of this
$\mathbb{P}^1$, it must hold at~$e_{u^\prime}$.

The second assertion is a special case of the contrapositive of the first.
\end{proof}

Our proof of Theorem~\ref{thm:poset} follows from the following local
isomorphism result, which explains the geometric significance of
interval pattern embeddings.

\begin{Theorem}
\label{thm:iso}
Let $\Phi$ be an interval pattern embedding of $[u,v]$ into
$[x,w]$. Then there exists a scheme-theoretic isomorphism ${\mathcal
N}_{u,v}\cong {\mathcal N}_{x, w}$.  Therefore, affine
neighborhoods of $X_v$ and $X_w$ respectively at $e_u$ and
$e_x$ are isomorphic up to cartesian products with affine
space.
\end{Theorem}

Note that if $\Phi$ does not induce an interval isomorphism of $[u,v]$
and $[x,w]$, then ${\mathcal N}_{u,v}$ and ${\mathcal N}_{x, w}$
are not isomorphic since their dimensions differ.

See Example~\ref{exa:proof_example} below for an illustration of the
arguments in the following proof.

\noindent
\emph{Proof of Theorem~\ref{thm:iso}:}
Let $I=\{i_1<\ldots < i_m\}$ be the embedding 
indices of $\Phi$ and
let $\{1,2,\ldots, n\}\setminus I =\{j_1<\ldots <j_{n-m}\}$. 

\begin{lemma}
\label{lemma:equals}
For $1\leq d\leq n-m$, 
\begin{equation}
\label{eqn:equals}
\# \{k: k\leq j_d \mbox{ and } \Phi(u)(k)\geq \Phi(u)(j_d)\}
=\# \{k: k\leq j_d \mbox{ and } w(k)\geq w(j_d)\}.
\end{equation}
\end{lemma}
\begin{proof}
Let $s_{\alpha_1}\cdots s_{\alpha_{\ell(v)-\ell(u)}}$ be a reduced expression
of $v^{-1}u$ and consider the corresponding product of (possibly
non-simple) transpositions $t_1\cdots t_{\ell(v)-\ell(u)}$ where
$t_{j}=(i_{\alpha_j} \leftrightarrow i_{\alpha_j +1})$.

Note that since $\Phi$ is an interval pattern embedding of $[u,v]$
into $[x,w]$, as we successively multiply $w$ on the right by the
transpositions $t_1, t_2,\ldots, t_{\ell(v)-\ell(u)}$, each
transposition drops the Coxeter length by exactly~1 (because each
decreases the Coxeter length, and the total drop in length is
$\ell(v)-\ell(u)$). Also observe that for any permutation $p$ and
transposition $t=(a \leftrightarrow b)$, $\ell(pt)=\ell(p)-1$ if and
only if $p(a)>p(b)$ and there does not exist an index $k$, $a<k<b$,
such that $p(a)>p(k)>p(b)$.

For each $c$, $0\leq c\leq\ell(v)-\ell(u)$, define the permutation
$w^{(c)}$ by $w^{(c)}:=wt_1\cdots t_c$.  We can now check that each
$w^{(c)}$ satisfies
$$\# \{k: k\leq j_d \mbox{ and } w^{(c)}(k)\geq w^{(c)}(j_d)\} =\#
\{k: k\leq j_d \mbox{ and } w^{(c-1)}(k)\geq w^{(c-1)}(j_d)\}$$ for
all $d$.  If this equation were to fail for any $d$, it would have
to be the case that $t_c=(a \leftrightarrow b)$ with $a<j_d<b$ and
$w^{(c-1)}(a)>w^{(c-1)}(j_d)>w^{(c-1)}(b)$.

Since $w^{(\ell(v)-\ell(u))}=\Phi(u)$, the lemma follows by induction.
\end{proof}

The following lemma shows the vanishing of certain coordinates
at all points of $N_{\Phi(u), w}$.

\begin{lemma}
\label{lemma:zeroed}
Let $g=(g_{ij}) \in {\mathcal N}_{\Phi(u), w}$.  Then for each 
$1\leq d\leq n-m$ we have $g_{\Phi(u)(j_d), j_d}=1$, 
$g_{a, j_d}=0$ for any $a\neq \Phi(u)(j_d)$ and
$g_{\Phi(u)(j_d), b}=0$ for any $b\neq j_d$.
\end{lemma}
\begin{proof}
Since $J_{\Phi(u)}$ vanishes on $g$, $g_{\Phi(u)(j_d), j_d}=1$,
$g_{a, j_d}=0$ for  $a< \Phi(u)(j_d)$ and
$g_{\Phi(u)(j_d), b}=0$ for $b>j_d$. It remains to check the cases
$a> \Phi(u)(j_d)$ and $b< j_d$.

Now we check the case $b<j_d$. Since $w(j_d)=\Phi(u)(j_d)$, both have
a ``1'' in position $(w(j_d), j_d)=(\Phi(u)(j_d), j_d)$. Moreover, by
Lemma~\ref{lemma:equals}, $r_{w(j_d), j_d}^{w}=r_{\Phi(u)(j_d),
j_d}^{\Phi(u)}$. Let this common integer be~$S$.  Then $r_{w(j_d), j_d
-1}^{w}=S-1$ since the ``$1$'' in position $(w(j_d),j_d)$ causes the
rank to increase by~$1$ as one moves from $(w(j_d), j_d -1)$ to
$(w(j_d), j_d)$.  Hence the $S\times S$ minors of the southwest
$(n-w(j_d)+1) \times (j_d-1)$ submatrix of $g$ vanishes.  Furthermore,
$r_{\Phi(u)(j_d),j_d-1}^{\Phi(u)}=S-1$, so there are $S-1$ rows in
this submatrix with an entry of ``1'' in their rightmost nonzero
columns; these rightmost columns are distinct from each other.  The
generic matrix $Z^{(\Phi(u))}$ has a ``0'' in row $\Phi(u)(j_d)$ in
these columns; if some other entry in row $\Phi(u)(j_d)$ has a nonzero
entry, then the submatrix would have $S$ linearly independent columns,
a contradiction.  Therefore, $g_{\Phi(u)(j_d),b}=0$ for any $b<j_d$.

The proof for the case $a>\Phi(u)(j_d)$ is similar.\end{proof}

Define the (algebraic) map $\Psi:{\mathcal N}_{\Phi(u), w}\to u\sigma
u^{-1}(\Omega^\circ_u)$ as the projection which deletes the columns
$j_1,\ldots, j_{n-m}$ and rows $w(j_1),\ldots, w(j_{n-m})$ from an
element $g\in {\mathcal N}_{\Phi(u), w}$. This map is well defined,
and, by Lemma~\ref{lemma:zeroed}, injective.

Next, we show that the image of $\Psi$ is actually inside ${\mathcal
N}_{u, v}$.  This amounts to verifying that $\Psi(g)$ satisfies the
southwest rank conditions corresponding to the minors generating
$I_{u, v}$.  Consider the southwest submatrix of $\Psi(g)$ with
northeast corner $(a,b)$; let $b^\prime=i_b$ and
$a^\prime=\Phi(u)(i_{u^{-1}(a)})$ be the corresponding indices for
$g$.  Observe that by the definition of the rank matrix,
$r^{\Phi(u)}_{a^\prime b^\prime}=r^u_{ab}+c_{ab}$ where $c_{ab}$
equals the number of positions of the form $(w(j_d), j_d)$ southwest
of $(a^\prime, b^\prime)$, which by Lemma~\ref{lemma:equals} is the
number of positions of the form $(\Phi(u)(j_d), j_d)$ weakly southwest
of $(a^\prime, b^\prime)$. But by Lemma~\ref{lemma:zeroed}, removing
row $\Phi(u)(j_d)$ and column $j_d$ must drop the rank by exactly 1,
since those rows and columns have a ``1'' at $(\Phi(u)(j_d),j_d)$ and
``0'' everywhere else.  Therefore, deleting the aforementioned rows
and columns precisely drops the rank at $(a^\prime,b^\prime)$ by
precisely $c_{ab}$, so $\Psi(g)$ satisfies exactly the rank conditions
for ${\mathcal N}_{u,v}$.  Therefore the image of $\Psi$ is inside
${\mathcal N}_{u,v}$.

On the other hand, given a point in ${\mathcal N}_{u,v}$ we can add
back these deleted rows and columns. This is clearly the inverse map
to $\Psi$, so it follows that $\Psi$ is an isomorphism from
$\mathcal{N}_{\Phi(u),w}$ to $\mathcal{N}_{u,v}$.\qed

We illustrate this theorem (and its proof) by the following example.

\begin{Example}
\label{exa:proof_example}
Let $v=35142, u=13524$ and $w=\underline{5} 8 9 \underline{7}
\underline{1} \underline{6} 2 \underline{3} 4$ be as in
Example~\ref{exa:35142var}, with $\Phi$ indicated by the underlined
positions, so $x=\Phi(u)==189573264$.  The intervals $[u,v]$ and
$[x,w]$ are easily checked to be isomorphic.  A generic matrix of
$u\sigma u^{-1}(\Omega^\circ_u)$ has the form
$$\left(
\begin{matrix}
1 & 0 & 0 & 0 & 0 \\
z_{41} & 0 & 0 & 1 & 0 \\
z_{31} & 1 & 0 & 0 & 0 \\
z_{21} & z_{22} & 0 & z_{24} & 1 \\
z_{11} & z_{12} & 1 & 0 & 0 
\end{matrix}
\right),
\mbox{ while }
v= 
\left(
\begin{matrix}
0 & 0 & 1 & 0 & 0 \\
0 & 0 & 0 & 0 & 1 \\
1 & 0 & 0 & 0 & 0 \\
0 & 0 & 0 & 1 & 0 \\
0 & 1 & 0 & 0 & 0 
\end{matrix}
\right)\!\!, \mbox{ and }
R^v = 
\left(
\begin{matrix}
1 & 2 & 3 & 4 & 5 \\
1 & 2 & 2 & 3 & 4 \\
1 & 2 & 2 & 3 & 3 \\
0 & 1 & 1 & 2 & 2 \\
0 & 1 & 1 & 1 & 1 
\end{matrix}
\right).$$
The reader can check that:
\[\mathcal{N}_{u,v}\cong \Spec(\frac{\mathbb{C}[z_{11}, z_{12}, z_{21}, z_{22}, z_{24}, z_{31}, z_{41}]}{\langle z_{11}, z_{21}, z_{22}, z_{41}\rangle})\cong {\mathbb A}^3.\]
Now,
$$
w=
\left(
\begin{matrix}
0 & 0 & 0 & 0 & 1 & 0 & 0 & 0 & 0  \\
0 & 0 & 0 & 0 & 0 & 0 & 1 & 0 & 0  \\
0 & 0 & 0 & 0 & 0 & 0 & 0 & 1 & 0  \\
0 & 0 & 0 & 0 & 0 & 0 & 0 & 0 & 1  \\
1 & 0 & 0 & 0 & 0 & 0 & 0 & 0 & 0  \\
0 & 0 & 0 & 0 & 0 & 1 & 0 & 0 & 0  \\
0 & 0 & 0 & 1 & 0 & 0 & 0 & 0 & 0  \\
0 & 1 & 0 & 0 & 0 & 0 & 0 & 0 & 0  \\
0 & 0 & 1 & 0 & 0 & 0 & 0 & 0 & 0  
\end{matrix}
\right)
\mbox{ and }
R^w=
\left(
\begin{matrix}
1 & 2 & 3 & 4 & 5 & 6 & 7 & 8 & 9  \\
1 & 2 & 3 & 4 & 4 & 5 & 6 & 7 & 8  \\
1 & 2 & 3 & 4 & 4 & 5 & 5 & 6 & 7  \\
1 & 2 & 3 & 4 & 4 & 5 & 5 & 5 & 6  \\
1 & 2 & 3 & 4 & 4 & 5 & 5 & 5 & 5  \\
0 & 1 & 2 & 3 & 3 & 4 & 4 & 4 & 4  \\
0 & 1 & 2 & 3 & 3 & 3 & 3 & 3 & 3  \\
0 & 1 & 2 & 2 & 2 & 2 & 2 & 2 & 2  \\
0 & 0 & 1 & 1 & 1 & 1 & 1 & 1 & 1
\end{matrix}
\right).
$$
A generic matrix in $x\sigma x^{-1}(\Omega^\circ_x)$ has the form
$$\left(
\begin{matrix}
1      & 0 &      0 & 0      & 0  & 0     & 0      & 0 & 0  \\
y_{81} & 0 &      0 & 0      & 0  & 0     & 1      & 0 & 0  \\
y_{71} & 0 &      0 & 0      & 0  & 1     & 0      & 0 & 0  \\
y_{61} & 0 &      0 & 0      & 0  & y_{66}& y_{67} & 0 & 1  \\
y_{51} & 0 &      0 & 1      & 0  & 0     & 0      & 0 & 0  \\
y_{41} & 0 &      0 & y_{44} & 0  & y_{46}& y_{47} & 1 & 0  \\
y_{31} & 0 &      0 & y_{34} & 1  & 0     & 0      & 0 & 0  \\
y_{21} & 1 &      0 & 0      & 0  & 0     & 0      & 0 & 0  \\
y_{11} & y_{12} & 1 & 0      & 0  & 0     & 0      & 0 & 0  
\end{matrix}
\right). \mbox{ Then }
\left(
\begin{matrix}
1 & 0 & 0 & 0 & 0 \\
y_{71} & 0 & 0 & 1 & 0\\
y_{51} & 1 & 0 & 0 & 0\\
y_{41} & y_{44} & 0 & y_{46} &1 \\
y_{31} & y_{34} & 1 & 0 & 0
\end{matrix}
\right)
$$
\noindent
is the result of applying $\Psi$ to this generic matrix, since $\Psi$
removes columns $2, 3, 7,$ and $9$ and rows $2, 4, 8,$ and $9$.  The
map $\Psi$ is injective since, by Lemma~\ref{lemma:zeroed}, $y_{11}$,
$y_{12}$, $y_{21}$, $y_{47}$, $y_{61}$, $y_{66}$, $y_{67}$, and
$y_{81}$ all equal $0$ for any point in $\mathcal{N}_{x,w}$.

Now examine the rank conditions imposed by $I_w$, which can be read
from the rank matrix $R^w$.  We have that $y_{31}=y_{41}=0$.
Furthermore, since $r^w_{6,5}=3$, all $4\times 4$ minors of the
southwest $4\times 5$ submatrix (of an element of $\mathcal{N}_{x,w}$)
vanish.  It follows that $y_{44}=y_{41}=0$.  The remaining rank
conditions imposing no further equations, it follows that ${\mathcal
N}_{x,w}\cong {\mathbb A}^3$, in agreement with
Theorem~\ref{thm:poset}. \qed
\end{Example}

\noindent
\emph{Conclusion of the proofs of Theorem~\ref{thm:poset} and
Corollary~\ref{cor:brpa}:} Let $[u,v]\in \mathfrak{S}$ be as in the
statement of the theorem.  Suppose $[u,v]\prec [x,w]$. We may assume
that this is a covering relation. There are two cases. If $\Phi$ is an
embedding of $v$ into $w$ and $x=\Phi(u)$ then the result follows from
Theorem~\ref{thm:iso}. In the other case, $[x,w]=[u^\prime,v]$ where
$u^\prime\leq u$, and the conclusion holds by
Lemma~\ref{lemma:semicont}.

The corollary is the contrapositive of the theorem, with
$A_\mathcal{P}$ being a generating set for~$\mathcal{I}_\mathcal{P}$.
\qed

\section{Computing with Kazhdan-Lusztig ideals}

Next we turn to computing properties of Kazhdan-Lusztig ideals
$I_{u,v}$ using {\tt Macaulay~2}. Various computations have been
automated in our scripts {\tt Schubsingular}. We illustrate the main
computations through specific examples here.  Our computations are over
$\mathbb{Q}$ and valid for any field of characteristic 0; similar
computations can be made in other characteristics.

These computations can be used to conjecture generators
$A_\mathcal{P}$ for the order ideal $\mathcal{I}_\mathcal{P}$ for
various particular properties $\mathcal{P}$.

To begin explaining the computations, let us re-examine
Example~\ref{exa:35142var}, in which $u=13254$ and $v=35142$.

\begin{verbatim}
R = QQ[z11, z12, z13, z14, z15, 
       z21, z22, z23, z24, z25,
       z31, z32, z33, z34, z35, 
       z41, z42, z43, z44, z45, 
       z51, z52, z53, z54, z55]; n=5;  -- graded diagonal term order on G 
                                       -- below
G = matrix({{  1,   0,   0,   0, 0},
            {z41,   0,   1,   0, 0},
            {z31,   1,   0,   0, 0},
            {z21, z22, z23,   0, 1},
            {z11, z12, z13,   1, 0}}); -- generic matrix in opposite 
                                       -- Schubert cell of u
Rank = matrix({{1, 2, 3, 4 ,5}, 
               {1, 2, 2, 3, 4}, 
               {1, 2, 2, 3, 3}, 
               {0, 1, 1, 2, 2}, 
               {0, 1, 1, 1, 1}});      -- rank matrix of v
Jlist = trim(sum(flatten(for i from 0 to n-1 list 
                          for j from 0 to n-1 list minors(Rank_(n-i-1,j)+1, 
                           submatrix(G, {(n-i-1)..n-1}, {0..j}))))) 
\end{verbatim}
The last line computes the generators for $I_{u,v}$. \qed

The following problem is the next step in all of our further
calculations, so an explicit combinatorial answer would speed up these
computations.
\begin{problem}
\label{prob:1}
Find a Gr\"{o}bner basis for $I_{u,v}$.
\end{problem}
It seems plausible to us that the defining determinants of $I_{u,v}$
are a Gr\"{o}bner basis with respect to any graded diagonal term order
in any characteristic.  (This has been verified computationally for
$n=4,5$, with respect to the canonical graded diagonal term order as
used above).  After re-writing our matrices upside down, this
conjecture would generalize the conclusion for matrix Schubert
varieties found in~\cite{grobGeom}. (In the process of preparing this
report, we mentioned this possibility to A.~Knutson, who informed us
that he had independently discovered this generalization.)


Finding a Gr\"{o}bner basis is the first step towards computing a free
resolution of $I_{u,v}$. We point out that explicitly finding a
(minimal) free resolution for even the special case of matrix Schubert
varieties is an open problem of considerable interest that has been
solved only in certain special cases.  We believe that free
resolutions for Kazhdan-Lusztig ideals in general, or even further
special cases, are also of interest.

When an ideal of a polynomial ring is homogeneous under some positive
grading, there is a unique minimal free resolution which is a
subcomplex of every free resolution.  (A positive grading is one for
which the degree $0$ piece is $\ourfield$ and the remainder of the
ring has positive degree.)  Unfortunately, $I_{u,v}$ is not in general
a homogeneous ideal under the naive grading with $\operatorname{deg} \
z_{ij} = 1$ for all $i,j$.  However, there is a positive grading of
$\ourfield[\zvect^{(u)}]$ under which $I_{u,v}$ is homogeneous.

\begin{lemma}
Let $e_1,\ldots,e_n$ be generators of the group $\mathbb{Z}^n$.  Under
the multi-grading where the variable $z_{n-u(i)+1,j}$ has degree
$e_i-e_j$, every variable in $\ourfield[\zvect^{(u)}]$ has a degree
which is a positive sum of the degrees $e_{i+1}-e_i$ for $1\leq i\leq
n-1$, and $I_{u,v}$ is homogeneous.
\end{lemma}

\begin{proof}
We can assign the degree $e_{u^{-1}(i)}-e_j$ to the entry in row $i$
and column $j$ of the generic matrix $Z^{(u)}$ to get this grading,
since a $0$ can be assigned any degree, the $1$'s, which are at
$(u(j),j)$ for $1\leq j\leq n$, are assigned the degrees $e_j-e_j=0$,
and the variables $z_{ij}$ are assigned the degrees specified.
Therefore, the minor of $Z^{(u)}$ using rows $i_1,\ldots,i_k$ and
columns $j_1,\ldots,j_k$ will be homogeneous of degree $\sum_{m=1}^k
e_{u^{-1}(i_m)} - \sum_{m=1}^k e_{j_m}$.  Since $I_{u,v}$ is generated
by minors of $Z^{(u)}$, this proves $I_{u,v}$ is homogeneous under
this grading.

For positivity, note that $Z^{(u)}$ has a variable $z_{n-u(i)+1,j}$ at
the position $(u(i),j)$ only if $j<i$.
\end{proof}

If we wish to have a $\mathbb{Z}$-grading, we can coarsen the above
grading by sending $e_{i+1}-e_i$ to $1$, which sends $e_i-e_j$ to
$i-j$; this coarser grading is also positive.  Also, note that
our isomorphism $\Psi$ of Theorem~\ref{thm:iso} is compatible with the map
of multigradings sending $e_{\phi_t}$ to $e_t$.  (Our multigrading
secretly comes from the action of $T$ on $G/B$.)

\begin{problem}
\label{prob:minfreeres}
Determine a (minimal) free resolution of $I_{u,v}$.
\end{problem}

The Betti number $\beta_i(I_{u,v})$ is the rank of the free module $F_i$ in
a minimal free resolution
$$0\leftarrow F_0 \leftarrow F_1 \leftarrow \ldots \leftarrow F_i
\leftarrow \ldots \leftarrow 0.$$ If we also keep track of the
(multi)-degrees of the generators of $F_i$, we get (multi-)graded
Betti numbers.  One problem which may be an easier step towards
finding a free resolution is the following.

\begin{problem}
\label{prob:3}
Give a combinatorial method for finding the (multi-graded) betti
numbers of $I_{u,v}$.
\end{problem}

It would then be left only to determine the maps.  Note that answers
to both of these problems are known to depend on the characteristic of
the field.

Continuing further with our re-examination of Example~\ref{exa:35142var}, let
us compute a minimal free resolution with respect to the aforementioned
grading on $I_{u,v}$.

\begin{verbatim}
S = QQ[z11, z12, z13,
       z21, z22, z23,
       z31,          
       z41, Degrees=>{{3},{2},{1},
                      {4},{3},{2},
                      {1},
                      {2}}];  -- S homogeneous with respect to 
                              -- the nonstandard grading

Jlist = substitute(Jlist, S)  -- convert Jlist into an ideal of S

Resl  = res(S^1/Jlist) -- free resolution as an S-module

       1      5      9      7      2
o2 = S  <-- S  <-- S  <-- S  <-- S  <-- 0

      0      1      2      3      4      5


o2 : ChainComplex

Resl.dd_1 -- minimal generators of Jlist

o3 = | z11 z12z31+z13z41 z21 z13z22-z12z23 z22z31+z23z41 |

             1       5
o3 : Matrix S  <--- S

Resl.dd_2 -- first syzygies, etc

o4 = {3} | 0    -z12z31-z13z41 0    -z21 -z13z22+z12z23 0             
     {3} | -z23 z11            -z22 0    0              -z21          
     {4} | 0    0              0    z11  0              z12z31+z13z41 
     {4} | -z31 0              z41  0    z11            0             
     {4} | z13  0              z12  0    0              0             

                       -z22z31-z23z41 0               0              |
                       0              0               0              |
                       0              -z13z22+z12z23  -z22z31-z23z41 |
                       0              z21             0              |
                       z11            0               z21            |
                                                  
              5       9
o4 : Matrix S  <--- S

betti(Resl) -- degrees are for the grading mentioned above

o5 = total: 1 5 9 7 2
         0: 1 . . . .
         1: . . . . .
         2: . 2 . . .
         3: . 3 1 . .
         4: . . 2 . .
         5: . . 4 1 .
         6: . . 2 2 .
         7: . . . 2 .
         8: . . . 2 1
         9: . . . . 1
\end{verbatim}
(In $\tt{Schubsingular}$, this resolution is obtained using
{\tt minresKL($\{2,4,0,3,1\}, \{0,2,1,4,3\}$)}. Take note of the ``computer
indexing'' which shifts the numbers down by~1.)

Although the above grading is natural in some ways, it also causes
some problems as our ring $\ourfield[\zvect^{(u)}]$ is not generated
in degree 1.  Therefore, we are also interested in the following
problem.
\begin{problem}
Find an explicit characterization of the pairs $(u,v)\in S_n \times S_n$
for which $I_{u,v}$ is homogeneous under the standard grading 
${\rm deg} \ z_{ij} =1$.
\end{problem}

An obvious subset of such pairs consists of those where the essential set
of $v$ is contained inside the ``staircase'' defined by the $1$s of $u$.

\section{Calculations for singularity invariants}

In this section, we discuss various semicontinuously stable properties
${\mathcal P}$. In many instances, we present or conjecture a
nonrecursive combinatorial description of ${\mathcal I}_{\mathcal P}$ or
$A_{\mathcal P}$.  In other instances, we explain
computational ({\tt Macaulay~2}) aspects of the effort to find them.

\subsection{Smoothness}

The problem of determining the singular locus of $X_w$ was solved
independently by \cite{billey.warrington,Cortez, KLR, manivel1}; see
also the earlier work~\cite{Gasharov}.  We can restate this problem in terms of
interval pattern avoidance as that of finding a full set of
generators for the ideal $\mathcal{I}_{\rm singular}$.  It is not
difficult to verify that the following is a restatement of the singular
locus theorem. (In what follows, 
we use the convention that the segment ``$j\cdots i$'' means $j, j-1, j-2,\ldots, i+1,i$. 
In particular, if $j<i$ then the segment is empty.)

\begin{Theorem}
The order ideal $\mathcal{I}_{\mathrm{singular}}$ is minimally
generated by the following families of intervals:
\begin{enumerate}
\item $\big[(a+1)a\cdots 1(a+b+2)\cdots(a+2),\ \ \ (a+b+2)(a+1)a\cdots 2 (a+b+1)\cdots(a+2) 1\big]$ for all integers $a,b>0$.
\item $\big[(a+1)\cdots 1 (a+3) (a+2) (a+b+4)\cdots(a+4),\ \ \ (a+3)(a+1)\cdots 2(a+b+4)1(a+b+3)\cdots(a+4)(a+2)\big]$ for all integers $a,b\geq0$.
\item $\big[1(a+3)\cdots2(a+4), \ \ \  (a+3)(a+4)(a+2)\cdots312 \big]$ for all integers $a>1$.
\end{enumerate}
\end{Theorem}

\begin{Example}
\label{exa:maxsing_461253}
We compute the singular locus of $X_{461253}$ by the above theorem.
After calculating the pairs that arise for $S_n$, $n\leq 6$, we have
only an embedding of $[1324,3412]$ into the first four positions of
$461253$, and two embeddings of $[13254, 35142]$, one excluding
position $3$ and the other excluding position $4$.  Therefore,
$X_{461253}$ is singular at $e_{u^\prime}$ if and only if
$u^\prime\leq 142653$, $u^\prime\leq 241365$, or $u^\prime\leq 143265$
in Bruhat order.  In other words, the singular locus decomposes as
\[{\rm sing}(X_{461253})=X_{142653}^{\circ}\cup X_{241365}^{\circ}\cup
X_{143265}^{\circ}.\]\qed
\end{Example}

The property ${\mathcal P}=$``singular'' has the special feature that
the set of permutations appearing as the top element of intervals in
${\mathcal I}_{{\rm singular}}$ is the order ideal generated by $4231$
and $3412$ in the partial order given by \emph{classical} pattern
avoidance, where ``$u$ is smaller than $v$'' if $u$ classically embeds
into $v$.  This is the theorem of Lakshmibai and Sandhya \cite{LS90}
stated in Section~1.

This special feature does not hold in general for all semicontinuously
stable properties.  (Compare Conjecture~\ref{conj:nonGorlocus} below
with Example~\ref{exa:Gor_not_pattern}).  On the other hand, Billey
and Braden~\cite{BilBrad} have given a geometric explanation of why
ordinary pattern avoidance characterizes smoothness.  One would like
to know which other semicontinuously stable properties ${\mathcal P}$
ordinary pattern avoidance characterizes.  We also wonder if some
feature of the combinatorics of ${\mathcal I}_{\mathcal P}$ might
characterize when ordinary pattern avoidance actually suffices.

\subsection{Kazhdan-Lusztig polynomials}

Associated to each pair of permutations 
$v,w\in S_n$ with $v\leq w$ is the
{\bf Kazhdan-Lusztig polynomial} $P_{v,w}(q)\in {\mathbb
N}[q]$. Although these polynomials have an elementary recursive
definition~\cite{KL}, their combinatorics is difficult to understand.
This is one motivation for studying the singularities of Schubert
varieties.  

Geometrically, $P_{v,w}(q)$ is the Poincar\'{e} polynomial for the
local intersection cohomology of $X_w$ at $e_v$.  (Given this is a
topological invariant, the problem of calculating Kazhdan-Lusztig
polynomials only makes sense over the field $\mathbb{C}$.)  In
particular (in type $A$), \linebreak $P_{v,w}(q)=1$ if and only if
$X_w$ is smooth at $e_v$.  Although dimensions of local intersection
cohomology groups are not in general (upper or lower) semicontinuous,
a result of Irving~\cite{Irving} (see also~\cite{braden.macpherson})
shows that they behave in an upper semicontinuous manner on Schubert
varieties.  Therefore, we can study Kazhdan-Lusztig polynomials using
interval pattern avoidance.

It is an longstanding, well-known open problem to find (hopefully
nonrecursive) manifestly positive, combinatorial rules for
$P_{v,w}(q)$. Such rules are only known in a limited cases,
essentially those where a semi-small resolution of singularities is
known~\cite{BW-321hex, Las-vex}.  The following corollary of
Theorem~\ref{thm:iso} appears to be new.  It generalizes a lemma of
P.~Polo~\cite[Lemma 2.6]{Polo} (see also \cite[Theorem 6]{BilBrad})
which states the case where the embedding $\Phi$ is in consecutive
positions or the entries in the positions of $\Phi$ are numerically
consecutive.

\begin{corollary}
\label{cor:KL_cor}
Suppose $[u,v]$ and $[x,w]$ are isomorphic because of an interval
pattern embedding.  Then $P_{x,w}(q)=P_{u,v}(q)$.
\end{corollary}

Lusztig's {\bf interval conjecture} asserts that
$P_{a,b}(q)=P_{v,w}(q)$ whenever the Bruhat order intervals $[a,b]$
and $[v,w]$ are isomorphic as posets; this conjecture is discussed
with further references in~\cite{brenti, bjorner.brenti}.  
Corollary~\ref{cor:KL_cor}
confirms a new (albeit, very special) case of the conjecture.

\begin{Example}
\label{exa:KL_diff}
Ordinary pattern avoidance does not suffice for
Corollary~\ref{cor:KL_cor}. Consider $v=4231$ embedding into
$\underline{52}3\underline{4} \underline{1}$ at the indicated
positions. Let $u=2143\leq v$. So $\Phi(u)=21354$. Then $P_{u,v}=1+q$
while $P_{\Phi(u),w}=1+2q+q^2$.  (Here we have made use of Goresky's
tables for Kazhdan-Lusztig polynomials~\cite{goresky}.)
\end{Example}

At present, we have no counterexample to the analogue of the interval
conjecture for any of the semicontinuously stable properties studied
in this paper.  Therefore, the much stronger assertion that
$\mathcal{N}_{a,b}\cong\mathcal{N}_{u,v}$ whenever $[a,b]$ and $[u,v]$
are isomorphic as intervals in Bruhat order remains possible, though
in our opinion extremely unlikely.

As with all of the numerical invariants studied here, understanding of
where they increase is a basic issue of interest. This suggests a new
incremental formulation of the aforementioned Kazhdan-Lusztig
polynomial problem.
\begin{problem} 
Let $\mathcal{P}_{k,\ell}$ to be the property ``the coefficient of
$q^\ell$ in $P_{u,v}(q)$ is at least $k$'' (or equivalently
``$\dim_\ourfield IH^\ell_{e_u}(X_v)\geq k$'').  Determine
$\mathcal{I}_{\mathcal{P}_{k,\ell}}$ for various values of $k$
and~$\ell$.
\end{problem}

As with all such problems in this paper, as a first step we can
formulate the computational challenge of analyzing ${\mathcal
I}_{{\mathcal P}}^{(n)}$ for small~$n$, where ${\mathcal I}_{{\mathcal
P}}^{(n)}$ is the set of intervals from ${\mathcal I}_{\mathcal P}$
from $S_m$ for $m\leq n$.

Unfortunately, there appears to be no known applicable method for
computing the ranks of local intersection cohomology groups directly
from the equations defining a variety.

\subsection{The Gorenstein property and Cohen-Macaulay type}

A local ring $(R,\mathfrak{m}, \Bbbk)$ is said to be {\bf Cohen-Macaulay}
if $\Ext_R^i(\Bbbk, R)=0$ for $i\leq\dim R$; it is {\bf
Gorenstein} if, in addition, $\dim_{\Bbbk}
\Ext_R^{\dim R}(\Bbbk, R)=1$.  A variety is Cohen-Macaulay
(respectively Gorenstein) if the local ring at every point is
Cohen-Macaulay (respectively Gorenstein).  Using the Kozsul complex on
a regular sequence, one can show that every regular local ring is
Gorenstein; hence smooth varieties are Gorenstein.  See
\cite{Bruns-Herzog} for details. 

In \cite{WY:Gor} we characterized the Schubert varieties which
are Gorenstein at all points.  Here is a reformulation of the main
result of that paper purely in terms of interval pattern
avoidance. 

\begin{Theorem}
\label{thm:Gor_char}
The Schubert variety $X_w$ is Gorenstein if and only if $w$
avoids the following intervals
\begin{enumerate}
\item $\big[(a+1)a\cdots 1(a+b+2)\cdots(a+2), \ \ \ (a+b+2)(a+1)a\cdots 2 (a+b+1)\cdots(a+2) 1\big]$ for all integers $a,b>0$ such that $a\neq b$.
\item $\big[(a+1)\cdots 1 (a+3) (a+2) (a+b+4)\cdots(a+4),\ \ \ (a+3)(a+1)\cdots 2(a+b+4)1(a+b+3)\cdots(a+4)(a+2)\big]$ for all integers $a,b\geq0$, with either $a>0$ or $b>0$.
\end{enumerate}
\end{Theorem}

As is explained in~\cite{WY:Gor}, the above theorem shows that a
Schubert variety is Gorenstein if and only if the generic points of
its singular locus are. This is closely related to the following
conjecture, which we now restate in the terminology of
interval pattern avoidance.

\begin{conjecture}
\label{conj:nonGorlocus}
The order ideal $\mathcal{I}_{\mathrm{not \ Gorenstein}}$ is generated by
the following families:
\begin{enumerate}
\item $\big[(a+1)a\cdots 1(a+b+2)\cdots(a+2), \ \ \ (a+b+2)(a+1)a\cdots 2 (a+b+1)\cdots(a+2) 1\big]$ for all integers $a,b>0$ such that $a\neq b$.
\item $\big[(a+1)\cdots 1 (a+3) (a+2) (a+b+4)\cdots(a+4),\ \ \ (a+3)(a+1)\cdots 2(a+b+4)1(a+b+3)\cdots(a+4)(a+2)\big]$ for all integers $a,b\geq0$, with either $a>0$ or $b>0$.
\end{enumerate}
\end{conjecture}

The components of the singular locus whose generic points are not
Gorenstein are clearly in the non-Gorenstein locus.  \emph{A priori},
it is possible for some Schubert variety to be non-Gorenstein at some
non-generic points outside of these components known to be
non-Gorenstein.  The content behind the conjecture is that this does
not occur.  Translating this geometric assertion into the
combinatorics of interval pattern avoidance, we get the
conjecture that the generators of $\mathcal{I}_{\mathrm{not \
Gorenstein}}$ is the subset of the minimal generators $(u,v)$ of
$\mathcal{I}_{\mathrm{singular}}$ for which $X_v$ is not Gorenstein at
$e_u$.  Using the description of neighborhoods of generic points of
the singular locus given in~\cite{Cortez,manivel2}, we arrived at the
above conjecture, which has been partially verified by computations as
described below.

The {\bf Cohen-Macaulay type} of a local Cohen-Macaulay ring is defined
to be \linebreak $\dim_{\Bbbk} \Ext_R^{\dim R}(\Bbbk, R)$. This is a numerical
refinement of the binary question of whether a variety is Gorenstein
at a point.  In our case, where the ring
$R=\ourfield[\zvect^{(u)}]/I_{u,v}$ is given as a quotient of a
polynomial ring $S=\ourfield[\zvect^{(u)}]$, Cohen-Macaulay type can be 
calculated as
the last non-zero Betti number by the following well-known argument.
We recall it here for completeness, as we could not find an explicit
reference.

We have that $\Ext_R^{\dim R}(\Bbbk,R)=\Ext_S^{\dim R}(\Bbbk, R)$.
Now let $K_\bullet$ be the Kozsul complex which is a free resolution
of $\Bbbk$ (over $S$); we then can calculate the $\Ext$ module using
the definition $\Ext_S^{\dim R}(\Bbbk, R)=H^{\dim
R}(\Hom_S(K_\bullet,R))$.  Since $K_\bullet$ is self-dual of length
$\dim S$, we have that $H^{\dim R}(\Hom_S(K_\bullet,R))=H_{\dim S-\dim
R}(K_\bullet\otimes R)=\Tor_{\codim R}^S(\Bbbk, R)$.  Now we can
calculate $\Tor_{\codim R}^S(\Bbbk, R)$ using a free resolution of $R$
as an $S$-module, and the Cohen-Macaulay type can be calculated as the
last (since Schubert varieties are Cohen-Macaulay) non-zero Betti
number of $R$ as an $S$-module.

Computation gives us the following partial check of our conjecture.
Note the computations have only been done in characteristic 0.

\begin{Proposition}
Conjecture~\ref{conj:nonGorlocus} is true for $n\leq 6$.
\end{Proposition}
\begin{proof}
The conjecture is vacuously true for $n\leq 4$ and for $n=5$ the result
follows from~\cite[Corollary~1]{WY:Gor},
since every non-Gorenstein Schubert variety for $n=5$ has only one
component in its singular locus, so every singular point is
non-Gorenstein.  The verification for $n=6$ is by computer.
\end{proof}

The computation verifying the conjecture for
$n=6$ is already somewhat involved, as the following
example shows.  

\begin{Example}
\label{exa:461253}
Let $v=461253$ as in the previous example.  The conjecture states that
$X_{461253}$ is not Gorenstein at $e_{u^\prime}$ if and only if
$u^\prime\leq 241365$ or $u^\prime\leq 143265$ in Bruhat order.  There
are 24 elements of $S_n$ in the interval between $v$ and $id$
(inclusive). Of those, 9 are in fact smaller that $241365$ or $143265$
in the Bruhat order, namely:
$$142365, 124365, 132465, 142356, 123465, 124356, 132456, \idperm.$$
All of the remainder are larger than $u=123546$ in the Bruhat
order. Therefore, one only needs to compute a free resolution for
$\mathcal{N}_{u,v}$.  {\tt Macaulay 2} reveals the following.
\begin{verbatim}
      1      7      21      35      35      21      7      1
o4 = S  <-- S  <-- S   <-- S   <-- S   <-- S   <-- S  <-- S  <-- 0

     0      1      2       3       4       5       6      7      8

o4 : ChainComplex
\end{verbatim}
The last nonzero free module in the resolution is rank 1, 
agreeing with the conjecture. (Using {\tt Schubsingular} the function {\tt locuscompute} automates
the Bruhat order analysis needed in the above
computation.)\qed
\end{Example}

The most general problem for this invariant is to calculate the
Cohen-Macaulay type of $X_v$ at $e_u$.  As with the problem of
calculating Kazhdan-Lusztig polynomials, this problem can also be
reformulated in an incremental form as follows.

\begin{problem}
\label{prob:cansheafrankgen}
Let $\mathcal{P}_k$ be the property ``canonical sheaf of $X_v$ has
rank at least $k$''.  Find the generators for the ideal
$\mathcal{I}_{\mathcal{P}_k}$ for all $k$.
\end{problem}

It may be particularly interesting to understand what changes in the
singularity structure cause the Cohen-Macaulay type to change,
leading to the following question.

\begin{problem}
\label{prob:typechange}
Characterize pairs $(u,v)$ such that the Cohen-Macaulay type of $X_v$
at $e_u$ is larger than the Cohen-Macaulay type $X_v$ at
$e_{u^\prime}$ for all $u^\prime$ with $u<u^\prime\leq v$.
\end{problem}

One important value for $u$ is the identity permutation, since, by
semicontinuity, the Cohen-Macaulay type of $X_v$ will be largest
at that point.  Therefore, the following case of
Problem~\ref{prob:cansheafrankgen} is of particular interest.

\begin{problem}
\label{prob:typeid}
Find the Cohen-Macaulay type of $X_v$ at $\idperm$.
\end{problem}

For $n\leq 4$, all $X_w$ are Gorenstein, so these types are all equal to~1.
For $n=5$, the four non-Gorenstein Schubert varieties 
\[X_{53241}, X_{35142}, X_{42513}, X_{52431}\]
all have type $2$ at the identity. Using, 
{\tt Schubsingular} we determined that $X_{624351}$ is the unique
Schubert variety having type $4$ at the identity, while
\[X_{361542}, X_{426153}, X_{623541}, X_{625431}, X_{532614}, X_{632451}, 
X_{643251}\] all have type $3$ there. (We use the command {\tt
cansheafrank($\{2,5,0,4,3,1\}, \{0,1,2,3,4,5\}$)} in {\tt
Schubsingular} to make the computation for $X_{361542}$ and similar
commands for the others.)  The remaining Schubert varieties have type
$1$ or $2$ at the identity. These can be distinguished by
Theorem~\ref{thm:Gor_char}.

\subsection{Factoriality and the Class Group}

A variety is said to be {\bf factorial} if the local ring at every
point is a unique factorization domain.  Since regular local rings are
unique factorization domains, every smooth variety is factorial.
Furthermore, all unique factorization domains are Gorenstein.

M.~Bousquet-M\'elou and S.~Butler have characterized factorial Schubert
varieties by the following theorem: 

\begin{Theorem}[\cite{butler.melou}]
The Schubert variety $X_w$ is factorial if and only if $w$ classically avoids
$4231$ and interval avoids $[3142,3412]$.
\end{Theorem}

The considerations that led to Conjecture~\ref{conj:nonGorlocus} also
lead to the following conjecture.

\begin{conjecture}
\label{conj:nonfactorialconj}
The order ideal $\mathcal{I}_{\mathrm{not \ factorial}}$ is generated by
the following families:
\begin{enumerate}
\item $\big[(a+1)a\cdots 1(a+b+2)\cdots(a+2), \ \ \ (a+b+2)(a+1)a\cdots 2 (a+b+1)\cdots(a+2) 1\big]$ for all integers $a,b>0$.
\item $\big[(a+1)\cdots 1 (a+3) (a+2) (a+b+4)\cdots(a+4), \ \ \ (a+3)(a+1)\cdots 2(a+b+4)1(a+b+3)\cdots(a+4)(a+2)\big]$ for all integers $a,b\geq0$.
\end{enumerate}
\end{conjecture}

As with the Gorenstein property, there is also an invariant which
measures how far a local ring is from being a unique factorization
domain, the (local Weil) class group.  A local ring is a unique
factorization domain if and only if the class group is trivial.  We do
not know of an algorithm to compute these class groups, or to
otherwise check Conjecture~\ref{conj:nonfactorialconj}.

\subsection{Multiplicity}

The {\bf multiplicity} of a local ring $(R,\mathfrak{m}, \Bbbk)$ is
the degree of the projective tangent cone
$\operatorname{Proj}(\operatorname{gr}_{\mathfrak{m}} R)$ as a
subvariety of the projective tangent space
$\operatorname{Proj}(\operatorname{Sym}^{\ast}\mathfrak{m}/\mathfrak{m}^2)$.
Equivalently, if the Hilbert-Samuel polynomial of $R$ is $a_n x^{n}+
a_{n-1} x^{n-1} +\cdots + a_0$, then the multiplicity of $R$ is
$n!a_n$.  Given a scheme $X$ and a point $p$, the multiplicity of $X$
at $p$, usually denoted $\operatorname{mult}_p(X)$, is the
multiplicity of the local ring
$(\mathcal{O}_{X_p},\mathfrak{m}_{p},\Bbbk)$.

It is an open problem to find a manifestly positive combinatorial rule
for the multiplicity of $X_v$ at the point $e_u$.  Several such rules
are known for the case where $v$ is a Grassmannian permutation
\cite{Kratt, KreiLak, LakWey, RosZel}.  The analogue of
Corollary~\ref{cor:KL_cor} for multiplicity is the following:

\begin{corollary}
\label{cor:mult_eq}
Suppose we have an interval pattern embedding of $[u,v]$ into $[x,w]$.
Then
$\operatorname{mult}_{e_x}(X_w)=\operatorname{mult}_{e_u}(X_v)$.
\end{corollary}

\begin{Example}
Ordinary pattern avoidance also fails for
Corollary~\ref{cor:mult_eq}. With the same choice of
$u,v,w$ and $\Phi$ as in Example~\ref{exa:KL_diff}, we have 
$\operatorname{mult}_{e_{\Phi(u)}}(X_w)=5$ while
$\operatorname{mult}_{e_u}(X_v)=2$.
\end{Example}

We reformulate the problem of finding a rule for multiplicity 
in an incremental form as follows.

\begin{problem}
\label{prob:multgen}
Let $\mathcal{P}_k$ be the property ``multiplicity of $X_v$ is at
least $k$''.  Find the generators for the ideal
$\mathcal{I}_{\mathcal{P}_k}$ for all $k$.
\end{problem}

Now we look at how multiplicity can be computed.  In our coordinates
for $\mathcal{N}_{u,v}$, $e_u$ is the point where $z_{ij}=0$ for all
$i,j$.  In this case where the local ring is the localization of a
ring $S/J$ at the maximal ideal $\mathfrak{m}$ given by the vanishing
of all the variables, the associated graded ring
$\operatorname{gr}_\mathfrak{m} S/J$ is isomorphic to the localization
of $S/J^\prime$ (at the maximal ideal given by the vanishing of all
variables), where $J^\prime=\langle f^\prime\mid f\in J\rangle$, with
$f^\prime$ defined to be the sum of all terms of minimal degree in
$f$.  To calculate the degree of $J^\prime$, we can then calculate a
Gr\"obner basis of $J^\prime$; this entire process of finding
$J^\prime$ and finding a Gr\"obner basis for it can be accomplished in
one step by finding the Gr\"obner basis (or an initial ideal) of $J$
with respect to a term order that chooses a {\it lowest} degree term.
Note that we are now using the grading where each variable has degree
1 rather than the grading discussed in section 4.
 
In {\tt Macaulay~2}, we can simulate a term order choosing a lowest
degree term by homogenizing the generators of $I_{u,v}$ using a new
variable $t$ (or by replacing the ``$1$''s by ``t''s in the matrix
corresponding to $u\sigma u^{-1}(\Omega_{u}^{\circ})$) and using a
term order that refines the partial order by degree in $t$. We can
then compute the initial ideal, send $t$ to $1$ and calculate the
degree of the resulting (monomial) ideal to find multiplicity.

\begin{Example}
Returning back to Example~\ref{exa:35142var}, let us calculate the multiplicity
of $X_{35142}$ at $e_{13254}$. 

\begin{verbatim}
i2 : St  = QQ[t, z11, z12, z13,
              z21, z22, z23,
              z31,
              z41, MonomialOrder=>Eliminate 1]; n=5;
i4 : Rank = matrix({{1, 2, 3, 4 ,5},
                    {1, 2, 2, 3, 4},
                    {1, 2, 2, 3, 3},
                    {0, 1, 1, 2, 2},
                    {0, 1, 1, 1, 1}}); -- rank matrix of w
i5 : Gt  = matrix({{  t,   0,   0,   0, 0},
                   {z41,   0,   t,   0, 0},
                   {z31,   t,   0,   0, 0},
                   {z21, z22, z23,   0, t},
                   {z11, z12, z13,   t, 0}});
i6 : Jlist = trim(sum(flatten(for i from 0 to n-1 list
                                for j from 0 to n-1 list
                                  minors(Rank_(n-i-1,j)+1,
                                  submatrix(Gt, {(n-i-1)..n-1}, {0..j})))))
i7 : GBlist = gb(Jlist);
i8 : LTlist = leadTerm(gens(GBlist)); -- gives in(J_{v,w})
i9 : S   = QQ[z11, z12, z13, 
              z21, z22, z23,
              z31,
              z41];
i10 : f = map(S, St, {1, z11, z12, z13, z21, z22, z23, z31, z41});
i11 : ELTlist = f(LTlist); -- gives in(J_{v,w}')
i12 : Dlist = degree(ideal(ELTlist))
o12 = 2
\end{verbatim}
Hence ${\rm mult}_{e_{13254}}(X_{35142})=2$. (In {\tt Schubsingular}, 
this calculation is automated by the command 
{\tt mult(\{2,4,0,3,1\},\{0,2,1,4,3\})}.)\qed
\end{Example}

One possible method for solving the problem of finding multiplicities
would be to find a combinatorial description for the initial ideals
resulting from the above algorithm, under a particularly good choice
of term order.  Under particularly good conditions, the set of pipe
dreams for the matrix Schubert variety $\pi^{-1}(X_v)$ counts the
multiplicity; it is shown by this method in \cite{Woo} that the
multiplicity of $X_{n23\cdots(n-1)1}$ at $e_{\idperm}$ is the Catalan
number $C_{n-2}$, and conjectured there that this is the highest
multiplicity (at any point) on any Schubert variety in
${\rm Flags}({\mathbb C}^n)$.

We also have the analogues of Problems~\ref{prob:typechange}
and~\ref{prob:typeid} for multiplicity.

\begin{problem}
Characterize pairs $(u,v)$ such that the multiplicity
of $X_v$ at $e_u$ is larger than the multiplicity
$X_v$ at $e_{u^\prime}$ for all $u^\prime$ with $u<u^\prime\leq v$.
\end{problem}

\begin{problem}
Find the multiplicity of $X_v$ at $\idperm$.
\end{problem}

Based on our calculations so far, the projection of the
ideal for the property ``multiplicity of $X_v$ at $e_u$ is at least
3'' onto the second factor is an order ideal in the partial order
given by ordinary pattern avoidance.  Both geometric and combinatorial
explanations of this phenomenon, if it indeed holds in general, would
be interesting.

\subsection{Final remarks and summary for $n=5$} 
In this report we have discussed several semicontinuously stable
invariants of Schubert varieties. We present a compact summary
for $n=5$, which can also be verified using {\tt Schubsingular}.

\begin{Proposition}
\label{prop:n=5}
\begin{itemize}
\item $X_{52341}$ is Gorenstein, has multiplicity 5 below 21354, multiplicity
1 where it is nonsingular, and multiplicity 2 everywhere else.

\item The 4 non-Gorenstein $X_w$ have multiplicity 3 and canonical sheaf
rank 2 where singular.

\item All other singular $X_w$ have multiplicity 2 where singular.
\end{itemize}
\end{Proposition}

{\tt Schubsingular} provides algorithms to compute such facts for
larger~$n$. However, already for $n=6$ the situation is complex enough 
that we refer the reader to the software.

     There are other interesting cases of invariants not considered here.
As one example, 
a local ring is a {\bf complete intersection} if it is the quotient of a
regular local ring by a regular sequence.  Remarkably, this is
actually a homological property which is independent of the ambient
ring.  In our case, $X_v$ is locally a complete intersection at $e_u$
if and only if the first Betti number for $\mathcal{N}_{u,v}$ is equal
to $\binom{n}{2}-\ell(v)$.

Not being locally a complete intersection is a semicontinuously stable property;
indeed, the difference between the first Betti number of
$\mathcal{N}_{u,v}$ and $\binom{n}{2}-\ell(v)$ is upper semicontinous.
We conclude with the following question, which was raised independently 
by B.~Hassett, R.~Joshua and B.~Sturmfels:
\begin{problem}
Which Schubert varieties $X_w$ are local complete intersections?
\end{problem}

\section*{Acknowledgements}
We thank M.~Haiman and F.~Sottile for (separately) suggesting that
Bruhat-restricted pattern avoidance should have a geometric
explanation, inspiring the present study.  In addition, we thank
S.~Billey, J.~Carrell, A.~Cortez, A.~Knutson, V.~Lakshmibai, I.~Lankham,
E.~Miller, V.~Reiner and G.~Smith for helpful discussions.  This paper
was prepared in part during the authors' residence at the 2005 AMS
Summer Research Institute on Algebraic Geometry in Seattle.  In
addition, this work was partially completed while AY was an NSERC
supported member of the Fields' Institute during 2005-2006, and while
an NSF supported visitor at the Mittag-Leffler Institute during Spring
2005.

\end{document}